\pdfoutput=1 

\documentclass{IEEEtran}
\usepackage{amssymb,empheq,amsthm}
\usepackage{graphicx, subfigure, color}
\graphicspath{{Figure/}}
\usepackage[pdftex,linkcolor=black,citecolor=black,pdfborder={0 0 0}]{hyperref}
\allowdisplaybreaks
\usepackage{enumitem,cases}

\newcommand{\sR}{\mathbb{R}}
\newcommand{\sC}{\mathbb{C}}
\newcommand{\sN}{\mathbb{N}}

\newcommand{\sZ}{\mathbb{Z}}

\newcommand{\mi}{\mathbf{i}}
\newcommand{\mRe}[1]{~\Re(#1)}
\newcommand{\mIm}[1]{~\Im(#1)}
\usepackage{mathrsfs}
\newcommand{\gG}{\mathscr{G}}
\newcommand{\gV}{\mathscr{V}}
\newcommand{\gE}{\mathscr{E}}
\newcommand{\gA}{\mathscr{A}}

\newcommand{\mE}{\mathcal{E}}
\newcommand{\mA}{\mathcal{A}}
\newcommand{\mB}{\mathcal{B}}

\newcommand{\mC}{\mathcal{C}}
\newcommand{\mD}{\mathcal{D}}

\newcommand{\mJ}{\mathcal{J}}
\newcommand{\mP}{\mathcal{P}}
\newcommand{\mX}{\mathcal{X}}

\newcommand{\mQ}{\mathcal{Q}}
\newcommand{\mS}{\mathcal{S}}
\newcommand{\mU}{\mathcal{U}}
\newcommand{\mZero}{\mathbf{0}}
\newcommand{\mOne}{\mathbf{1}}
\newcommand{\mO}{\mathrm{O}}

\usepackage[f]{esvect}

\newcommand{\vq}{\mathbf{q}}

\newcommand{\diag}[1]{\mathrm{diag}\{#1\}}

\newcommand{\abs}[1]{|#1|}
\newcommand{\abss}[1]{\big|{#1}\big|}
\newcommand{\Abs}[1]{\bigg|{#1}\bigg|}
\newcommand{\norm}[1]{\|#1\|}

\newcommand{\norminf}[1]{\|#1\|_{\infty}}

\newcommand{\spec}[1]{\rho(#1)}

\usepackage{mathtools,extarrows}
\DeclarePairedDelimiter{\ceil}{\lceil}{\rceil}
\DeclarePairedDelimiter{\floor}{\lfloor}{\rfloor}

\newcommand{\ignore}[1]{}

\newcommand{\tand}{\text{~and~}}

\newcommand{\tfor}{\text{~for~}}
\newcommand{\tif}{\text{~if~}}
\newcommand{\tin}{\text{~in~}}

\usepackage{color,xcolor}
\definecolor{szu}{cmyk}{0.10,0.97,0.37,0.43}
\definecolor{blue2}{RGB}{0,112,192}
\definecolor{black2}{RGB}{128,128,128}

\def\smath#1{\text{\scalebox{.8}{$#1$}}}

\def\sfrac#1#2{\smath{\frac{#1}{#2}}}

\usepackage{tikz}

\theoremstyle{definition}
 
 \newtheorem{lem}{Lemma}
 
 \newtheorem{thm}{Theorem}
 \newtheorem{cor}{Corollary}
 
 \newtheorem{rem}{Remark}

\begin{document}
\title{Synchronization of Networked Harmonic Oscillators via Quantized Sampled Velocity Feedback}
\author{Jingyi Wang, Jianwen Feng, Yijun Lou, Guanrong Chen,~\IEEEmembership{Fellow,~IEEE}
\thanks{This work was supported in parts by the National Natural Science Foundation of China under Grants 61603260, 61873171 and 61872429.}
\thanks{J. Wang and J. Feng  are with the College of Mathematics and Statistics, Shenzhen University, Shenzhen 518060, China 
	(e-mail: wangjingyi@szu.edu.cn, fengjw@szu.edu.cn)}
\thanks{Y. Lou is with the Department of Applied Mathematics, The Hong Kong Polytechnic University,
 Hong Kong SAR, China (e-mail: yijun.lou@polyu.edu.hk)}
\thanks{G. Chen is with the Department of Electronic Engineering, City University of Hong Kong, Hong Kong SAR, China (e-mail: eegchen@cityu.edu.hk)}
}

\maketitle
 \begin{abstract}
In this technical note, we propose a practicable quantized sampled velocity data coupling protocol for synchronization of a set of harmonic oscillators.
The coupling protocol is designed in a quantized way via interconnecting the velocities encoded by a uniform quantizer with a zooming parameter in either a fixed or an adjustable form over a directed communication network.
We establish sufficient conditions for the networked harmonic oscillators to converge to a bounded neighborhood of the synchronized orbits with a fixed zooming parameter. 
We ensure the oscillators to achieve synchronization by designing the quantized coupling protocol with an adjustable zooming parameter.
Finally, we show two numerical examples to illustrate the effectiveness of the proposed coupling protocol.
\end{abstract}

\section{Introduction}\label{sec:introduction}

%
\IEEEPARstart{S}{ynchronization} phenomena are common in nature and society.
Understanding, describing and controlling synchronization have been an active research field in various academic disciplines, as surveyed by \cite{phyrep2002boccaletti.chaosync,ieeeproc2007olfati.consensus,auto2014dorfler:sync}.
In general, synchronization is a process in which the state of network-interconnected subsystems
converge to the same orbit driven by a designed coupling or control protocol \cite{phyrep2002boccaletti.chaosync}.

%
%
Synchronization of networked harmonic oscillators provides a basic model for studying the dynamics and control problems of complex dynamical networks, 
with significant practical applications, such as mobile robots \cite{auto2008ren:harm} and electrical networks \cite{auto2017Tuna}.
In the past decade, some effective coupling protocols have been presented from different perspectives.
For instance,
the oscillators with interactions in a continuous-time setting over fixed or switching network topologies were investigated in \cite{auto2008ren:harm},
while the discrete-time setting was studied and applied to synchronization control of multiple mobile robots in \cite{iet2010ballard:harm},
and their instantaneous interactions under fixed or switching topologies with presence or absence of leaders were considered in \cite{automatica2012zhou.pco}.
A distributed protocol was proposed in an impulsive form by using the relative position  information between the oscillator and its neighbors in \cite{ieeetac2018zhang:cho}. 
Synchronization can also be reached by directly utilizing delayed position states in \cite{ieeec2016song:cho}.
Recently, the synchronization problem was solved even by using noisy sampled-data in \cite{wang:auto2016nho,wang:auto2018nho}.
Nonlinear diffusive coupling can also achieve synchronization \cite{ieeetac2017liu:cho}.

%
%
On the other hand, almost all of the control systems are implemented digitally today, from large computer systems to small embedded processors,
for which sampling and quantization are fundamental tools.
In particular, in the study of communications and coupling protocols of oscillators, 
 transmitting data can be obtained via continuous, periodic or aperiodic sampling, which could be then sent and received in a quantized form.
With this motivation, stabilization of linear systems via quantized control was studied in the continuous-time \cite{tac2000brockett:qfs}, discrete-time \cite{auto2003liberzon} and switching \cite{ieeetac2017wakaiki:swls} settings, which may be subject to  external disturbances \cite{ieeetac2012sharon:cq}.
Moreover,  quantization techniques can be used to deal with the  stabilization of systems with limited measurement information \cite{ieeetac2001elia:slsli,ieeetac2003liberzon:slsli}.
The consensus of fixed \cite{scl2013xu:qm} or switched \cite{ijss2017zhu:qc} networks of multi-agent systems was studied based on quantized relative state information, 
as well as quantization of the absolute state information \cite{auto2011ceragioli:qc}.

To the best of our knowledge, few studies have been conducted on the synchronization problem of harmonic oscillators over a directed communication network via quantized feedback coupling by using only sampled velocity data.
%
%
The objective of this technical note is to present some novel coupling protocols using quantized sampled control to achieve synchronization of a network of harmonic oscillators.
First, a quantized coupling protocol with a fixed zooming parameter is proposed, which can guarantee the boundedness of the synchronization errors at periodic discrete-time instants.
Then, its modified version with an adjustable zooming parameter is designed to achieve complete synchronization for any initial values of the oscillators.


\subsubsection*{Notations}
Let $\mOne_n=[1,1,\ldots,1]^{\top} \in \sR^n$,
$\mZero_n=[0,0,\ldots,0]^{\top} \in \sR^n$,
$\mO_n = [0] \in \sR^{n \times n}$,
$\mO_{m,n} = [0] \in \sR^{m \times n}$ and $I_n$ be the $n-$dimensional identity matrix.
Use superscript $\top$ ($\ast$) to denote the transpose (conjugate transpose) of a vector or matrix.
For $x \in \sR$, the ceiling function $\ceil{ x }=\min\{k \in \sZ: k \geq x\}$ is the smallest integer not less than $x$,
and the floor function $\floor{ x}=\max\{k \in \sZ: k \leq x\}$ is the largest integer not greater than $x$.
Let $\mi=\sqrt{-1}$ denote the imaginary unit.
For a complex number $x$, $\abs{x}$ denotes its modulus, and $\mRe{x}$ and $\mIm{x}$ represent its real part and imaginary
part, respectively.
For a vector $x \in \sC^n$, the maximum norm ($l_{\infty}$-norm)  is denoted by $\norm{x}_{\infty}= \max_i \abs{x_i}$,
while the $l_2$-norm is denoted by $\norm{x}_2=\sqrt{x^{\ast}x}$.
For a square matrix $X \in \sC^{n \times n}$,
the spectrum of $X$ is denoted by $\sigma(X)$, that is, the set of eigenvalues of $X$,
while the spectral radius of $X$ is denoted by $\spec{X}=\max \{\abs{\lambda}: \lambda \in \sigma(X)\}$.
For a matrix $X \in \sC^{m \times n}$,
the maximum row-sum matrix norm, induced by the $l_{\infty}$-norm, is defined by
$\norm{X}_{\infty}=\max_i \sum_{j=1}^n\abs{x_{ij}}$, while the spectral norm, induced by the $l_2$-norm, is defined by $\norm{X}_2=\max\{\sqrt{\lambda}:\lambda \in \sigma (X^{\ast}X)\}$, that is, the largest singular value of $X$.

\section{Preliminaries and Problem  statement}\label{sec:preliminaries}

\subsection{Preliminaries}

To investigate the synchronization behavior, the oscillators are interconnected over a communication network described by a directed graph without loops. 
Let $\gG = (\gV, \gE, \gA)$ be a directed and connected graph without loops for
a set of oscillators $\gV=\{1,2, \ldots, n\}$, a set of edges $\gE\subseteq \gV \times \gV$
and the adjacency matrix $\gA = [a_{ij}] \in \sR^{n \times n}$, in which
$a_{ij}$ is a positive weight for all $(i,j) \in \gE$
if and only if oscillator $i$ can access or receive the information from oscillator $j$
and $a_{ij} = 0$ for all $(i,j) \notin \gE$.
Define the Laplacian matrix $L = [l_{ij}] \in \sR^{n \times n}$ associated with $\gG$
by $l_{ij} = -a_{ij}$ for $i \neq j$ and $l_{ii}=\sum_{i=1,j \neq i}^n a_{ij}$ for $i = j$.
A sequence of edges $(i_1,i_2)$, $(i_2,i_3)$, $\ldots$, with $i_j \in \gV$, is a directed path in $\gG$.
A directed graph is strongly connected if and only if any two distinct nodes of the graph can be connected via a directed path.

\begin{lem}\label{lem:LxiP}\cite{ieeetac2005ren:cons}
 If a directed graph $\gG$ is strongly connected, for which the associated Laplacian matrix is $L$, then
\begin{enumerate}
 \item $\mathbf{1}_n$ is a right eigenvector of $L$ associated with the eigenvalue $\lambda_1 = 0$ of multiplicity $1$, and all the other right eigenvalues $\lambda_2,\ldots,\lambda_n$ have positive real parts;
 \item if $\xi = [\xi_1, \xi_2, \ldots , \xi_n]^{\top}$ is a left eigenvector of $L$ associated with the eigenvalue $0$
(i.e. $\xi^{\top}L = 0$), then { $\xi_i > 0$} for all $i = 1, 2,\ldots, n$, and $\xi$ has multiplicity $1$;
 \item there exists a nonsingular matrix $P$, in which the first column is $\mathbf{1}_n$ and the first row of $P^{-1}$ is $\xi^{\top}$, such that $L = PJP^{-1}$ is the Jordan decomposition of $L$, where $J = \diag{0, \hat{J}}$
 and $\hat{J}$ is the Jordan upper diagonal block matrix corresponding to the nonzero eigenvalues $\lambda_r$ ($r = 2,\ldots, n$) of the matrix $L$.
\end{enumerate}
\end{lem}

For the above mentioned claim, assume that $\sum_{i=1}^n \xi_i = 1$, and denote $\bar{\xi}= \max \xi_i$,
$P = [p_1,p_2,\ldots,p_{n}] \in \sC^{n \times n}$, where $p_i \in \sC^{n \times 1}$, and
$P^{-1} = [\bar{p}_1^{\top},\bar{p}_2^{\top},\ldots,\bar{p}_{n}^{\top}]^{\top} \in \sC^{n \times n}$, where $\bar{p}_i \in \sC^{1 \times n}$.
Obviously, $p_1 = \mOne_n$, $\bar{p}_1 = \xi^{\top}$, $\bar{p}_ip_i = 1$ and $\bar{p}_ip_j = 0$ for all $i \neq j$.
Let $\Xi=p_1\bar{p}_1 \in \sR^{n \times n}$, $\hat{P} = [p_2,\ldots,p_{n}] \in \sC^{n \times n-1}$,
$\hat{P}^{\dagger} = [\bar{p}_2^{\top},\ldots,\bar{p}_{n}^{\top}]^{\top} \in \sC^{n -1 \times n}$,
so that $\hat{P}^{\dagger} \hat{P} = I_{n-1}$ and $\hat{P}^{\dagger} L \hat{P} = \hat{J}$,
where $\hat{J} = \diag{\hat{J}_{n_1},\hat{J}_{n_2},\ldots,\hat{J}_{n_r}} \in \sC^{(n-1) \times (n -1)}$,
in which $\hat{\mJ}_{n_k}$ is an $n_k \times n_k$ Jordan block corresponding to eigenvalue $\lambda_k$
%
with  geometric multiplicity $n_k$ and $n_1+n_2 +\cdots+n_r = n-1$.
Denote $\mP = \diag{P,P}$, 
$\mJ = \diag{J,J}$, $\hat{\mJ} = \diag{\hat{J},\hat{J}}$, $\hat{\mP} =\begin{bmatrix}
\hat{P} & \mO_{n,n-1} \\ \mO_{n,n-1} & \hat{P}
\end{bmatrix}$ and $\hat{\mP}^{\dagger} =\begin{bmatrix}
\hat{P}^{\dagger} & \mO_{n-1,n} \\ \mO_{n-1,n} & \hat{P}^{\dagger}
\end{bmatrix}$.


\begin{lem}\cite{book:ma}\label{lem:ma5610}
For given $A \in \sC^{n \times n}$ and $\varepsilon > 0$, there exists a matrix norm $\norm{\cdot}$ such
that $\rho(A) \leq \norm{A} \leq \rho(A) + \varepsilon$.
\end{lem}

\begin{lem}\cite{book:ma}\label{lem:ma5612}
Let $A \in \sC^{n \times n}$. Then, $\lim\limits_{k \rightarrow \infty} A^k = \mO_n$ if and only if $\rho(A) < 1$.
\end{lem}

\begin{lem}\label{lem:compstab}\cite{book1993parks:st}
Given a complex-coefficient polynomial of order two,
$g(s) = s^2 +(a+ b \mi) s + (c+d \mi)$,
where $a$, $b$, $c$, and $d$ are real constants. Then, $g(s)$ is stable, that is, all roots of $g(s)=0$ have negative real
parts,
 if and only if $a > 0$ and $abd + a^2c -d^2 > 0$.
\end{lem}



\subsection{Quantizer}


In this technical note,  the class of quantizers proposed in \cite{tac2000brockett:qfs,auto2003liberzon} is adopted.
Let $\mQ$ be a finite subset of $\sR$. A quantizer is a piecewise constant
function $q : \sR \rightarrow \mQ$. This implies geometrically that $\sR$ is divided
into a finite number of quantized regions $\{y \in \sR : q(y) = y_i, y_i \in \mQ\}$. For a quantizer $q$, there exist positive numbers $M$ and $\Delta$
with $M > \Delta$ such that
\begin{enumerate}[label={(\roman*)}]
 \item if $|y| \leq M$, then $|q(y) - y| \leq \Delta$;\label{cond:quani}
 \item if $|y|> M$, then $|q(y)|> M - \Delta$.\label{cond:quanii}
\end{enumerate}
Condition \ref{cond:quani} gives an upper bound for the quantization error when the quantizer does not saturate,
and condition \ref{cond:quanii} is used for  detecting  quantizer saturation.
Here, $M$ and $\Delta$ are referred to as the range of $q$ and the quantization error, respectively.
To achieve complete synchronization of a network of oscillators, a quantizer can be designed with a suitable parameter $\mu > 0$, with $q(y)= q_{\mu} (y) = \mu q (y/\mu)$.
The parameter $\mu$ is regarded as a zooming variable, and the measurement capability and precision accuracy of the quantizer can be adjusted by the zoom-in and zoom-out operations.

\subsection{Model description}

Consider $n$ harmonic oscillators with the control input in the following form:
\begin{align}\label{equ:hamui}
\left\{
\begin{aligned}
 \dot{r}_i(t) =&v_i(t),\\
 \dot{v}_i(t) =&- \omega^2 r_i(t) +u_i(t),
\end{aligned}\right.
\end{align}
where $r_i(t)$, $v_i(t)$ and $u_i(t) \in \sR$ are the position, velocity, and control input of oscillator $i$, respectively, $\omega$ is a positive gain,
and the network topology is to be specified.

The  objective is to design a coupling protocol for each
oscillator such that the networked oscillators can achieve synchronization in the sense that
$\lim\limits_{t\rightarrow\infty}||r_i(t)-r_j(t)||=0$ and
$\lim\limits_{t\rightarrow\infty}||v_i(t)-v_j(t)|| =0$
for any initial values $r_i(0),v_i(0) \in \sR$ with a desired  norm $\norm{\cdot}$.
In this note, the following coupling protocol, 
 using the quantization of the absolute sampled velocity data, is designed
\begin{multline}\label{equ:quantcont}
u_i(t) =- \sum_{j=1}^n a_{ij}(q_{\mu}(v_i(t_k)) - q_{\mu}(v_j(t_k))),\\ t\in[t_k,t_{k+1}),~i=1,2,\ldots,n,
\end{multline}
where $v_j(t_k)$ is the sampled velocity of  oscillator $j$, obtained at the sampling instants $t_k$, $k = 0,1,2,\ldots$, satisfying
$t_k = k \tau$, where $\tau>0$ is a fixed sampling period
and $q_{\mu}$ is a quantizer.
In this note, both fixed and adjustable zooming parameters are investigated for quantized control, as follows:
\begin{description}
 \item[Case 1] $q_{\mu}(v_i(t_k)):=\mu q\big(\frac{v_i(t_k)}{\mu}\big)$, where $\mu$ is a fixed zooming parameter;
 \item[Case 2] $q_{\mu}(v_i(t_k)):=q_{\mu(t_k)}(v_i(t_k))=\mu(t_k) q\big(\frac{v_i(t_k)}{\mu(t_k)}\big)$, where $\mu(t_k)$ is an adjustable zooming parameter.
\end{description}

\begin{rem}
 The orbit of a simple harmonic oscillator, when the control input is not applied, is an ellipse and the energy function $V(t)=r_i^2(t)+v_i^2(t)/\omega^2$  is a constant for $t \geq 0$ ($V'(t)=0$).
 If $\norm{[r_i(0), v_i(0)]^{\top}}_{\infty} < M_0$,
 one has $\norm{[r_i(t), v_i(t)]^{\top}}_{\infty} \leq \sqrt{1+\bar{\omega}^2}M_0 $ for $t \geq 0$, where $\bar{\omega}= \max\{\omega,{1}/{{\omega}}\}$ and $M_0$ is a constant.
\end{rem}


\section{Main Results}\label{sec:mainresult}

Let
$r(t)=[r_1(t), r_2(t),$ $\ldots,$ $r_n(t)]^{\top}$, $v(t)=[v_1(t),v_2(t),\ldots,v_n(t)]^{\top}$,
$X(t) =[ {r}(t)^{\top}, {v}(t)^{\top}]^{\top}$, 
$\vq_{\mu}(r(t))=[q_{\mu}({v}_1(t)),q_{\mu}({v}_2(t)),\ldots,q_{\mu}({v}_n(t))]^{\top}$,
$\mB=\diag{\mO_n,-L}$,
\begin{equation*}
\mA = \begin{bmatrix}
 \mO_n&I_n \\
 -\omega^2 I_n & \mO_n
 \end{bmatrix},
\mC(t_k)=
 \begin{bmatrix}
 \mZero_n \\
 -L \big(\vq_{\mu}({v}(t_k))-v(t_k)\big)
 \end{bmatrix},
\end{equation*}
$E = \exp(\mA \tau)+\int^{\tau}_{0}\exp(\mA s) \mB ds$, and
$F = \int^{\tau}_{0}\exp(\mA s) ds$.

Let $\mE = \mP^{-1}E\mP$ and $\hat{\mE} = \hat{\mP}^{\dagger}E\hat{\mP}$. 
By Lemma \ref{lem:ma5610}, one can construct a matrix norm $\norm{\cdot}_{\epsilon}~$ on $\sC^{(2n-2) \times (2n-2)}$ induced by a vector norm $\norm{\cdot}_{\epsilon}$ on $\sC^{2n-2}$, such that $\norm{\hat{\mE}}_{\epsilon} \leq \spec{\hat{\mE}} + \epsilon< 1$, as follows:
\begin{equation}\label{equ:norm}
\left\{\begin{aligned}
&\text{the vector norm~} \norm{\cdot}_{\epsilon}=\norm{\mD_{\epsilon}\mU \cdot}_{\infty},\\
&\text{the matrix norm~} \norm{\cdot}_{\epsilon} = \norm{\mD_{\epsilon} \mU \cdot \mU^{-1} \mD_{\epsilon}^{-1}}_{\infty},
\end{aligned}\right.
\end{equation}
where $\mU \in \sC^{(2n-2) \times (2n-2)}$ is a nonsingular matrix such that
$\mU^{-1}\hat{\mE}\mU = \diag{\bar{\mJ}_1,\bar{\mJ}_2,\ldots,\bar{\mJ}_r}$ is a Jordan matrix,
in which $\bar{\mJ}_k \in \sC^{n_k \times n_k }$ are Jordan blocks with $n_1+n_2+\cdots+n_r=2n-2$,
and $\mD_{\epsilon} = \diag{{\mD}_{\epsilon,1},{\mD}_{\epsilon,2},\ldots,{\mD}_{\epsilon,r}}$
in which ${\mD}_{\epsilon,k} = \diag{1,1/\epsilon,1/\epsilon^2,\ldots,1/\epsilon^{n_k-1}}$.

\subsection{Quantized feedback coupling with a fixed zooming parameter}
In this subsection,  consider the convergence of system \eqref{equ:hamui} under control \eqref{equ:quantcont} with a fixed zooming parameter $\mu$.
\begin{thm}\label{thm1}
Assume that the directed graph $\gG$ is strongly connected, with an arbitrary small $\varepsilon > 0$
 and  a large enough $M$ compared to $\Delta$ such that
 \begin{equation*}
 M > \sfrac{\norm{\hat{\mP} \mU^{-1}\mD_{\epsilon}^{-1}}_{\infty} \norm{\mD_{\epsilon}\mU\hat{\mP}^{\dagger}F\mB}_{\infty} }{ (1- \bar{\xi})(1 - \norm{\hat{\mE}}_{\epsilon})}\Delta,
 \end{equation*}
and that the sampling period satisfies 
\begin{equation}\label{them:tau}
	\tau \in \{\tau:\cot({\omega\tau}/{2}) > \phi_i,i=2,\ldots,n\},
\end{equation}
where
\begin{multline}\label{cond:tau}
\phi_i = \sfrac{\mRe{\lambda_i}\mIm{\lambda_i}^2+\mRe{\lambda_i}^3}{2 {\omega}\mRe{\lambda_i}^2} \\
+\sfrac{\sqrt{\big(\mRe{\lambda_i}\mIm{\lambda_i}^2+\mRe{\lambda_i}^3\big)^2 + 4 {{\omega^2}}\mRe{\lambda_i}^2 \mIm{\lambda_i}^2}}
 {2 {\omega}\mRe{\lambda_i}^2}.
\end{multline}
Then, the solutions $(r(t),v(t))$ of system \eqref{equ:hamui} under control \eqref{equ:quantcont} start from $(r(0),v(0))$ inside the set $\mS_1(\mu)$ will enter into the set $\mS_2(\mu)$ in finite time
\begin{align}\label{equ:T}
T = \ceil*{\log_{\norm{\hat{\mE}}_{\epsilon}}
	\sfrac{ \norm{\hat{\mP}\mU^{-1}\mD_{\epsilon}^{-1}}_{\infty} \norm{\mD_{\epsilon}\mU\hat{\mP}^{\dagger}F\mB}_{\infty}\Delta\varepsilon}{( 1- \bar{\xi} )(1-\norm{\hat{\mE}}_{\epsilon}) M}}\tau,
\end{align}
 where
 \begin{align}
 \mS_1(\mu)=&\big\{(r,v):\norm{\hat{\mP}^{\dagger}([r^{\top},v^{\top}]^{\top} -  [\gamma,\nu]^{\top} \otimes \mOne_n)}_{\epsilon}\notag\\
 &\leq \sfrac{ (1- \bar{\xi}) \mu M }{\norm{\hat{\mP}\mU^{-1}\mD_{\epsilon}^{-1}}_{\infty} } \big\},\label{equ:s1}\\
 \mS_2(\mu)=&\big\{(r,v):\norm{\hat{\mP}^{\dagger}([r^{\top},v^{\top}]^{\top} -  [\gamma,\nu]^{\top} \otimes \mOne_n)}_{\epsilon}\notag\\
     &\leq \sfrac{\norm{\mD_{\epsilon}\mU\hat{\mP}^{\dagger}F\mB}_{\infty}\Delta\mu(1+\varepsilon)}{1- \norm{\hat{\mE}}_{\epsilon}}\big\},\label{equ:s2}\\
 \gamma(t)=&\cos({\omega}t)\xi^{\top}r(0)+\sin({\omega}t)\xi^{\top}v(0)/{\omega},\label{equ:gamma}\\
 \nu(t)=&-\omega\sin({\omega}t)\xi^{\top}r(0)+\cos({\omega}t)\xi^{\top} v(0),\label{equ:nu}
 \end{align}
%
%
%
$\bar{\omega}= \max\{\omega,{1}/{{\omega}}\}$,
$\lambda_2,\lambda_3,\ldots,\lambda_n$ are non-zero eigenvalues of the Laplacian matrix $L$ of $\gG$,
$\xi=[\xi_1,\xi_2,\ldots,\xi_n]^{\top}$ is a left eigenvector of $L$ associated with the zero eigenvalue $\lambda_1$.
\end{thm}

\begin{IEEEproof}
By the properties of the quantizer $q$ and the Laplacian matrix $L$ of $\gG$, system \eqref{equ:hamui} with control \eqref{equ:quantcont}  together can be written as
\begin{align}\label{equ:sys1}
 \dot{X}(t) = \mA X(t) + \mB X(t_k) + \mC(t_k),\quad t\in[t_k,t_{k+1}).
\end{align}
And, note that $[\gamma(t),\nu(t)]^{\top}$ in \eqref{equ:gamma}
and \eqref{equ:nu} is the solution of the following equation:
\begin{equation}\label{equ:gammanu}
	\left\{
\begin{aligned}
\dot{\gamma}(t) =& \nu(t),\\
\dot{\nu}(t) =& -\omega^2 \gamma(t)
\end{aligned}	
	\right.
\end{equation}
 with initial value $[\xi^{\top}r(0),\xi^{\top}v(0)]^{\top}$. Simple computation yields that $\gamma^2(t)+\nu^2(t)/\omega^2=(\xi^{\top}r(0))^2 + (\xi^{\top}v(0))^2/\omega^2 = M_0^2/\omega^2$.
Obviously, $\mB X (t_k)$ and $\mC(t_k)$ are constant vectors. For $t \in [t_{k},t_{k+1})$,
integrating both sides of equation \eqref{equ:sys1} from $t_k$ to $t$, one obtains
\begin{align*}
 X (t)= &E(t,t_k) X (t_k) +F(t,t_k)\mC(t_k),
\end{align*}
where $E(t,t_k) = \exp(\mA (t-t_k))+\int^{t}_{t_k}\exp(\mA(t-s)) \mB ds$ and
$F(t,t_k) = \int^{t}_{t_k}\exp(\mA(t-s)) ds$.
Furthermore, for $t = t_{k+1}$, one has
\begin{equation}\label{equ:xtk}
 X (t_{k+1})= E(t_{k+1},t_k) X (t_k) + F(t_{k+1},t_k) \mC(t_k).
\end{equation}
Notice that $t_{k+1} - t_k = \tau$. So, $E(t_{k+1},t_k) =E$ and $F(t_{k+1},t_k) = F$ for all $k=0,1,2,\ldots$.

Next, in order to implement quantization, the condition $\norm{X(t_{k})}_{\infty} <\mu M$ must be satisfied for all $k=1,2,\ldots$.
First,
when $t=t_1$, one gets
\begin{align*}
 X (t_{1})= E X (t_0) + F\mC(t_0).
\end{align*}
Moreover, letting $\mX(t) = \mP^{-1} X (t) = [\mX_1(t),\ldots,\mX_{2n}(t)]^{\top} \in \sR^{2n}$ gives
\begin{equation}\label{equ:quantY0}
\mX (t_{1})= \mE\mX(t_0) + \mP^{-1} F \mC(t_0),
\end{equation}
where $\mE = \mP^{-1}E\mP$.
Denote $\bar{\mX}(t)=[\mX_{1}(t),\mX_{n+1}(t)]^{\top} \in \sR^2$ and $\hat{\mX}(t)=[\mX_{2}(t),\ldots,\mX_n(t),\mX_{n+2}(t),\ldots,\mX_{2n}(t)]^{\top} \in \sR^{2n-2}$.
By using the property of $P^{-1}$ (the first row of $P^{-1}$ is $\xi^{\top}$), 
one gets $\bar{p}_1 L (q_{\mu}({v}(t_0))-v(t_0))=0$,
so that equation \eqref{equ:quantY0} can be written as
\begin{subequations}\label{equ:mXt0}
	\begin{empheq}[left={\empheqlbrace}]{align}
		\bar{\mX}(t_1) =& [\gamma(t_1),\nu(t_1)]^{\top},\label{equ:mXt01}\\
		\hat{\mX}(t_1) =&\hat{\mE}\hat{\mX}(t_0) + \hat{\mP}^{\dagger} F \mC(t_0).\label{equ:mXt02}
	\end{empheq}
\end{subequations}
By the definition of $\mP$ ($p_1 = \mOne_n$), one obtains
\begin{align*}
X (t_1) = \mP \mX(t_1) = [\gamma(t_1),\nu(t_1)]^{\top} \otimes \mOne_n + \hat{\mP} \hat{\mX}(t_1).
\end{align*} 
By using $\norm{X(t_0)}_{\infty} \leq M \mu$, $M_0 \leq \bar{\xi}\mu M $ and the definition of $\mS_1(\mu)$,
it follows that
\begin{align*}
\norminf{X(t_1)} =& \norminf{[\gamma(t_1),\nu(t_1)]^{\top} \otimes \mOne_n + \hat{\mP} \hat{\mX}(t_1)}\\
\leq & M_0 + \norm{\hat{\mP}\mU^{-1}\mD_{\epsilon}^{-1} \mD_{\epsilon}\mU \hat{\mX}(t_1)}_{\infty}\\
\leq & \bar{\xi} \mu M + \norm{\hat{\mP}\mU^{-1}\mD_{\epsilon}^{-1}}_{\infty}\norm{\hat{\mX}(t_0)}_{\epsilon}\\
\leq & \mu M.
\end{align*}
Meanwhile, based on the definition of $\mS_2(\mu)$, from equation \eqref{equ:mXt02} and condition \eqref{them:tau}, it follows that
%
\begin{align*}
   & \norm{\hat{\mX}(t_1)}_{\epsilon} -\norm{\hat{\mX}(t_0)}_{\epsilon} \notag\\
\leq &(\norm{\hat{\mE}}_{\epsilon} - 1) \norm{\hat{\mX}(t_0)}_{\epsilon} + \norm{\hat{\mP}^{\dagger} F \mC(t_0)}_{\epsilon}\notag\\
< &0.
\end{align*}
%
%
%
In the same way, by mathematical induction and the definition of $\mS_1(\mu)$ and $\mS_2(\mu)$, one can show that when $$\frac{ (1- \bar{\xi}) \mu M }{\norm{\hat{\mP}\mU^{-1}\mD_{\epsilon}^{-1}}_{\infty} } \geq \norm{\hat{\mX} (t_{k})}_{\epsilon} \geq \frac{\norm{\mD_{\epsilon}\mU\hat{\mP}^{\dagger}F\mB}_{\infty}\Delta\mu}{1- \norm{\hat{\mE}}_{\epsilon}},$$ the inequality $\norm{X(t_{k})}_{\infty} <\mu M$ holds and the sequence $\{\norm{\hat{\mX}(t_k)}_{\epsilon}\}$ decreases.

Furthermore,  equation \eqref{equ:xtk} implies
\begin{multline*}
 X (t_{k})= E^k X (t_0) + E^{k-1} F \mC(t_0) + E^{k-2} F \mC(t_1)
 \\+ \cdots +E F \mC(t_{k-2}) + F\mC(t_{k-1}).
\end{multline*}
and
\begin{multline}\label{equ:quantY}
 \mX (t_{k})= \mE^k\mX(t_0) + \mE^{k-1} \mP^{-1} F \mC(t_0) + \mE^{k-2} \mP^{-1} F \mC(t_1)\\
 + \cdots+\mE\mP^{-1} F \mC(t_{k-2}) + \mP^{-1} F\mC(t_{k-1}).
\end{multline}
By using the property of $P^{-1}$ (the first row of $P^{-1}$ is $\xi^{\top}$), 
one can get $\bar{p}_1 L (q_{\mu}({v}(t_k))-v(t_k))=0$,
so that equation \eqref{equ:quantY} can be written as
\begin{equation}\label{equ:mX2}
\left\{\begin{aligned}
 \bar{\mX}(t_k) =& [\gamma(t_k),\nu(t_k)]^{\top},\\
 \hat{\mX}(t_k) =&\hat{\mE}^k\hat{\mX}(t_0) + \hat{\mE}^{k-1} \hat{\mP}^{\dagger} F \mC(t_0) + \hat{\mE}^{k-2} \hat{\mP}^{\dagger} F \mC(t_1)\\
 				 &+ \cdots+\hat{\mE}\hat{\mP}^{\dagger} F \mC(t_{k-2}) + \hat{\mP}^{\dagger} F\mC(t_{k-1}),
\end{aligned}\right.
\end{equation}
where $[\gamma(t),\nu(t)]^{\top}$ is defined by equation \eqref{equ:gammanu}.
From Lemma \ref{lem:ma5612}, Theorems 5.4.10 and 5.6.15 in \cite{book:ma}, and equation \eqref{equ:mX2},
one can prove that the vector sequence $\{\hat{\mX}(t_k)\}$ converges with respect to any norm if $\spec{\hat{\mE}} < 1$.

Hereinafter, we consider the case of $\spec{\hat{\mE}} < 1$.
{In this scenario, there exists a sufficiently small $\epsilon>0$ such that $\spec{\hat{\mE}} + \epsilon< 1$.}
Using the properties of the quantizer $\abs{q_{\mu}(v_i(t))-v_i(t)} \leq \Delta\mu$,
one has $\norm{\hat{\mP}^{\dagger}F \mC(t_k)}_{\epsilon} \leq \norm{\mD_{\epsilon}\mU\hat{\mP}^{\dagger}F\mB}_{\infty}\Delta\mu$
for $k = 1,2,\ldots$, and 
\begin{align}\label{equ:hatmxtk}
 \norm{\hat{\mX} (t_{k})}_{\epsilon}& \leq \norm{\hat{\mE}}_{\epsilon}^k\norm{\hat{\mX}(t_0)}_{\epsilon} + \norm{\hat{\mE}}_{\epsilon}^{k-1} \norm{\hat{\mP}^{\dagger}F \mC(t_0)}_{\epsilon}\notag\\
 &\quad+ \norm{\hat{\mE}}_{\epsilon}^{k-2} \norm{\hat{\mP}^{\dagger}F \mC(t_1)}_{\epsilon}
 + \cdots\notag\\
 &\quad+\norm{\hat{\mE}}_{\epsilon}\norm{\hat{\mP}^{\dagger} F \mC(t_{k-2})}_{\epsilon} + \norm{\hat{\mP}^{\dagger} F\mC(t_{k-1})}_{\epsilon}\notag\\
 & \leq \norm{\hat{\mE}}_{\epsilon}^k\norm{\hat{\mX}(t_0)}_{\epsilon}
 + \sfrac{\norm{\mD_{\epsilon}\mU\hat{\mP}^{\dagger}F\mB}_{\infty} \Delta\mu(1- \norm{\hat{\mE}}_{\epsilon}^k)}{1- \norm{\hat{\mE}}_{\epsilon}}\notag\\
 & \leq \norm{\hat{\mE}}_{\epsilon}^k\norm{\hat{\mX}(t_0)}_{\epsilon}
 + \sfrac{\norm{\mD_{\epsilon}\mU\hat{\mP}^{\dagger}F\mB}_{\infty} \Delta\mu}{1- \norm{\hat{\mE}}_{\epsilon}}.
\end{align}
By  $\norm{\hat{\mE}}_{\epsilon}< 1$, we further obtain
\begin{align}\label{equ:xtkbound}
 \limsup_{k \rightarrow \infty }\norm{\hat{\mX} (t_{k})}_{\epsilon}
 \leq \sfrac{\norm{\mD_{\epsilon}\mU\hat{\mP}^{\dagger}F\mB}_{\infty}\Delta\mu}{1- \norm{\hat{\mE}}_{\epsilon}}.
\end{align}
Using the definitions of $\mS_1(\mu)$ and $\mS_2(\mu)$, and inequality \eqref{equ:hatmxtk},
the solutions $r(t)$ and $v(t)$ of system \eqref{equ:hamui} under control \eqref{equ:quantcont} go to $\mS_2(\mu)$ from $\mS_1(\mu)$ by time $T$, which was defined in \eqref{equ:T}.

{Finally,  sufficient conditions are derived to ensure that $\spec{\hat{\mE}} < 1$.}
The characteristic polynomial of $\hat{\mE}$ is $p(x) = \Pi_{i=2}^{n} p_i(x)$,
where
\begin{align*}p_i(x)=\det
\begin{bmatrix}
 \cos(\omega\tau) - x & \frac{1}{{\omega}}\sin(\omega\tau) + \frac{\lambda_i }{\omega^2}(\cos(\omega\tau)-1) \\
 -{\omega}\sin(\omega\tau) & \cos(\omega\tau)-\frac{ \lambda_i }{{\omega}}\sin(\omega\tau) -x
\end{bmatrix}.
\end{align*}
So, $\spec{\hat{\mE}}<1$ if and only if $\abs{x} < 1$, where $x$ is the root of $p_i(x) = 0$ for all $i= 2,3,\ldots,n$.
By simply calculations, $p_i(x) =0$ can be rewritten as
\begin{multline}\label{equ:pix}
 p_i(x) = x^2 + \big(\sin(\omega\tau)\lambda_i/{\omega} - 2\cos(\omega\tau)\big) x \\
 + \big( 1 - \sin(\omega\tau)\lambda_i/{\omega}\big) = 0.
\end{multline}
If $\sin(\omega\tau) = 0$, then  $\abs{x} = 1$, so $\tau \neq \frac{k\pi}{{\omega}}$ for all $k\in \sN$, that is, $\sin(\omega\tau) \neq 0$ and $\cos(\omega\tau) = \pm 1$.
Let $x = (s+1)/(s-1)$. Then, equation \eqref{equ:pix} is transformed into
%
\begin{equation}\label{equ:pis}
s^2 + \sfrac{\lambda_i}{\omega} \cot({\omega\tau}/{2})s +
\cot^2 ({\omega\tau}/{2})-\sfrac{\lambda_i}{\omega}\cot({\omega\tau}/{2})= 0.
\end{equation}
It is easy to see from the property of the bilinear transformation  that $\abs{x}<1$ in equation \eqref{equ:pix} holds if and only if $\mRe{s}<0$ in equation \eqref{equ:pis} holds.
By Lemma \ref{lem:compstab}, one can get that $\mRe{s} < 0$ if and only if $$\frac{1}{{\omega}} \cot({\omega\tau}/{2})\mRe{\lambda_i} > 0$$ and
\begin{multline}\label{equ:tauinequ}
 {\omega} \mRe{\lambda_i}^2\cot^2({\omega\tau}/{2})
 - \left(\mRe{\lambda_i}\mIm{\lambda_i}^2 + \mRe{\lambda_i}^3\right)\\
 \times\cot({\omega\tau}/{2}) - {\omega} \mIm{\lambda_i}^2>0.
\end{multline}
Recall that $\mRe{\lambda_i} > 0$ for all $i=2,3,\ldots,n$,
and inequality \eqref{equ:tauinequ} holds by equation \eqref{cond:tau}.
So, if condition \eqref{cond:tau} holds, then $\spec{\hat{\mE}}<1$.
This completes the proof of Theorem \ref{thm1}.
 \end{IEEEproof}

\begin{rem}
 From Lemma \ref{lem:ma5612}, the matrix sequence $\{\hat{\mE}^k\}$ converges to $\mO_{2n-2}$ if and only if $\spec{\hat{\mE}} < 1$,
  and from  Lemma \ref{lem:ma5610}, one can see that $\spec{\hat{\mE}} < 1$ is a necessary condition for $\norm{\hat{\mE}} \leq \spec{\hat{\mE}} + \epsilon < 1$ to guarantee the convergence of the series $\sum_{k=0}^{\infty}\hat{\mE}^{k} \mP^{-1} F \mC(t_i)$.
  So, the case of $\spec{\hat{\mE}} < 1$ is only needed to be considered in the proof of Theorem \ref{thm1}.
 At the same time, the condition $\spec{\hat{\mE}} < 1$ is a necessary and sufficient condition to achieve synchronization for system \eqref{equ:hamui} under control \eqref{equ:quantcont}, in which $q(v(t_k))=v(t_k)$.
 This can be proved in the same way as the proof of Theorem \ref{thm1}, so it is omitted.
\end{rem}

\begin{rem}
Obviously, for any $S \in \sC^{(2n-2) \times (2n-2)}$,
one has
$\max_{\norm{x}_{\epsilon} =1}\norm{S x}_{\epsilon}
= \max_{\norm{\mD_{\epsilon}\mU x}_{\infty} =1}\norm{\mD_{\epsilon}\mU S x}_{\infty}
=\max_{\norm{D_{\epsilon}\mU x}_{\infty} =1}\norm{\mD_{\epsilon}\mU S \mU^{-1} D_{\epsilon}^{-1}D_{\epsilon}\mU x}_{\infty}
= \norm{S}_{\epsilon}$,
so the matrix norm $\norm{\cdot}_{\epsilon}$ can be induced by the vector norm $\norm{\cdot}_{\epsilon}$ in Theorem \ref{thm1}.
For a finite-dimensional real or complex vector space, all norms are equivalent.
Consequently, the convergence of a sequence of vectors in a finite-dimensional space is independent of the norm (see, e.g., Corollary 5.4.6. in \cite{book:ma}).
\end{rem}

\begin{rem}
In Theorem \ref{thm1}, the set $\mS_2$ is not an invariant region of system \eqref{equ:hamui}, but $(r(t),v(t)) \in \mS_2$, where $t= k\tau$ for all $k \geq T/\tau$ and $k \in \sN$.
By some elementary computations, one can verify that, for $t\in (t_k,t_{k+1})$,
$ X (t) = \mP \mX(t) = [\gamma(t),\nu(t)]^{\top} \otimes \mOne_n + \hat{\mP} \hat{\mX}(t)$,
$\hat{\mX} (t)=\hat{\mE}(t-t_k) \hat{\mX}(t_k) + \hat{\mP}^{\dagger} F(t,t_k) \mC(t_{k})$
and
$\norm{[r(t)^{\top},v(t)^{\top}]^{\top} -[\gamma(t),\nu(t)]^{\top} \otimes \mOne_n }_{\infty} \leq \norm{\hat{\mE}(t-t_k)}_{\epsilon} \norm{\hat{\mE}}_{\epsilon}^k\norm{\hat{\mX}(t_0)}_{\epsilon}
 + {\norm{\hat{\mE}(t-t_k)}_{\epsilon}\norm{\mD_{\epsilon}\mU\hat{\mP}^{\dagger}F\mB}_{\infty} \Delta\mu(1- \norm{\hat{\mE}}_{\epsilon}^k)}/{(1- \norm{\hat{\mE}}_{\epsilon})}
 +\norm{\mD_{\epsilon}\mU\hat{\mP}^{\dagger}F(t,t_k)\mB}_{\infty}\Delta\mu$. Obviously, the norms $\norm{\hat{\mE}(t-t_k)}_{\epsilon}$ and $\norm{F(t,t_k)}_{\epsilon}$ are time-varying but bounded,
 and the set $\mS_2$ is more important and useful than the invariant region of system \eqref{equ:hamui}.
 From equation \eqref{equ:xtkbound}, we can show that the norm
 $\norminf{[r(t_k)^{\top},v(t_k)^{\top}]^{\top} - \mOne_n \otimes [\gamma(t_k),\nu(t_k)]^{\top}}$ is bounded as $k \rightarrow \infty$,
 that is, the states of the synchronized oscillators converge to a bounded region of the  orbits.
\end{rem}

If the interconnected network is undirected, one has the following corollary.

\begin{cor}\label{cor:exzhou}
Assume that the graph $\gG$ is undirected and connected with an arbitrary $\varepsilon > 0$,
 and  large enough $M$ compared to $\Delta$ such that
 $$ M > \sfrac{\norm{\hat{\mP} \mU^{-1}\mD_{\epsilon}^{-1}}_{\infty} \norm{\mD_{\epsilon}\mU\hat{\mP}^{\dagger}F\mB}_{\infty} }{ (1- \bar{\xi})(1 - \norm{\hat{\mE}}_{\epsilon})}\Delta,$$
 for $\tau \in \{\tau:\cot({\omega\tau}/{2}) > {\lambda_i}/{\omega},i=2,3,\ldots,n\}$.
Then, the solutions $(r(t),v(t))$ of system \eqref{equ:hamui} under control \eqref{equ:quantcont} starting form $(r(0),v(0))$  inside the set $\mS_1(\mu)$ will enter into the set $\mS_2(\mu)$ in finite time.
 \end{cor}

\subsection{Quantized feedback coupling with an adjustable zooming parameter}
In this subsection, we consider the convergence of equation \eqref{equ:hamui} with an adjustable zooming parameter $\mu$ at sampling instants.

\begin{thm}\label{thm2}
 Assume that the directed graph $\gG$ is strongly connected with an arbitrarily small $\varepsilon > 0$,
 and that $M$ is large enough compared with $\Delta$ such that
 \begin{equation*}\label{equ:thm2:cond1}
  M > \max\Big\{2\Delta,\sfrac{\norm{\hat{\mP} \mU^{-1}\mD_{\epsilon}^{-1}}_{\infty} \norm{\mD_{\epsilon}\mU\hat{\mP}^{\dagger}F\mB}_{\infty} }{ (1- \bar{\xi})(1 - \norm{\hat{\mE}}_{\epsilon})}\Delta\Big\}.
 \end{equation*}
If $\cot({\omega\tau}/{2}) > \phi_i$ for all $i=2,\ldots,n$, where $\phi_i$ is defined by equation \eqref{cond:tau}.
 Then, there exists a right-continuous and piecewise-constant function $\mu(t)$ such that the solutions $[r_i(t),v_i(t)]^{\top}$ of system \eqref{equ:hamui} under control \eqref{equ:quantcont} exponentially converge to $[\gamma(t), \nu(t)]^{\top}$ in equations \eqref{equ:gamma} and \eqref{equ:nu}.
\end{thm}

\begin{IEEEproof}
To construct the adjustable zooming parameter $\mu$, we divide the proof into two steps.

 \textit{Step 1. The zooming-out stage,} in which one can increase $\mu$ to obtain a larger quantization range
 such that the quantizer can capture the output.

 Set the control law $u_i(t)=0$. Let $\mu(t_0) = 0$ and $\mu(t) = {k} \Delta$ for $t \in [t_k,t_{k+1})$.
 Then, there exists a ${k_0} \in \sN$ such that
 \begin{align}\label{equ:determinek}
 \Abs{\sfrac{v_i(t_{k_0})}{\mu(t_{k_0})}} \leq M - 2 \Delta, \tfor i=1,2,\ldots,n,
 \end{align}
 where $t_{k_0}=k_0 \tau$.
 This implies, by condition \ref{cond:quani}, that
\begin{equation*}
\left\{\begin{aligned}
 &\abss{q\big(\sfrac{v_i(t_{k_0})}{\mu(t_{k_0})}\big)} \leq M - \Delta,\\
 &\abss{q_{\mu(t_{k_0})}({v_i(t_{k_0})})- \Delta\mu(t_{k_0})} \leq M\mu(t_{k_0}).
\end{aligned}\right.
\end{equation*}
 In view of conditions \ref{cond:quani} and \ref{cond:quanii}, one has $(r(t_{k_0}),v(t_{k_0})) \in \mS_1({\mu(t_{k_0})})$.

 \textit{Step 2: The zooming-in stage,} in which one can decrease $\mu$ to obtain a smaller quantization error
 such that the solution $(r_i(t),v_i(t))$ of the system converges to $(\gamma(t),\nu(t))$.

 Let $u_i(t)$ be as in equation \eqref{equ:quantcont} in Case 2,
where $\mu(t) = \mu(t_{k_0}) = {k_0} \Delta $ for $t= [t_{k_0},t_{k_0}+T)$, and $T$ is given by equation \eqref{equ:T}.
By using Theorem \ref{thm1}, we have $(r(t_{k_0}+T),v(t_{k_0}+T)) \in \mS_2({\mu(t_{k_0})})$.

 For $t= [t_{k_0}+T,t_{k_0}+2T)$, let $\mu(t) = \mu(t_{k_0}+T) = \theta \mu(t_{k_0}) = \theta {k_0} \Delta $, where
 \begin{align*}
 \theta =\sfrac{\norm{\hat{\mP} \mU^{-1}\mD_{\epsilon}^{-1}}_{\infty} \norm{\mD_{\epsilon}\mU\hat{\mP}^{\dagger}F\mB}_{\infty} (1+\varepsilon)\Delta}{ (1- \bar{\xi})(1 - \norm{\hat{\mE}}_{\epsilon}) M}.
 \end{align*}
 Obviously, $\theta < 1$ by equation \eqref{equ:thm2:cond1}, and $\mS_1(\mu(t_{k_0}+T)) = \mS_2(\mu(t_{k_0}))$ by equations \eqref{equ:s1} and \eqref{equ:s2}.
 From Theorem \ref{thm1}, one has
 $(r(t_{k_0}+2T),v(t_{k_0}+2T)) \in \mS_2(\mu(t_{k_0}+T))$.

 By mathematical induction, 
for $t= [t_{k_0}+kT,t_{k_0}+(k+1)T)$, letting $\mu(t) = \mu(t_{k_0}+kT) = \theta^k \mu(t_{k_0}) = \theta^k {k_0} \Delta $, we obtain
$(r(t_{k_0}+(k+1)T),v(t_{k_0}+(k+1)T)) \in \mS_2(\mu(t_{k_0}+kT))$ and $\mS_2(\mu(t_{k_0}+kT)) = \mS_1(\mu(t_{k_0}+(k+1)T))$.

Based on the above analysis, the boundedness of $\norm{\hat{\mE}(t-t_k)}_{\epsilon}$ and $\norm{F(t,t_k)}_{\epsilon}$ for $t\in [t_k,t_{k+1})$,
imply that $\mu(t) \rightarrow 0$ and
$\norm{ [r(t)^{\top},v(t)^{\top}]^{\top} - \mOne_n \otimes [\gamma(t),\nu(t)]^{\top}}_{\infty} \rightarrow 0$ as $t \rightarrow \infty$.
This completes the proof of Theorem \ref{thm2}.
\end{IEEEproof}

\begin{rem}
By the proof of Theorem \ref{thm2}, one obtains the piecewise control law $u_i(t)$ and the zooming parameter $\mu(t)$ as follows:
\begin{equation}\label{equ:uimuconc}
\left\{\begin{aligned}
 &u_i(t) = 0 \tand \mu(t) = \ceil{{t}/{\tau}} \Delta, & \tif t <k_0 \tau,\\
 &u_i(t) \tin \eqref{equ:quantcont} \tand \mu(t) = \theta^{\floor{\frac{t-\tau k_0}{T}}} {k_0} \Delta, & \tif t \geq {k_0} \tau,
\end{aligned}\right.
\end{equation}
where $k_0$ is dependent on the initial values and  determined by equation \eqref{equ:determinek}.
\end{rem}

\section{Numerical Simulations}\label{sec:numericalsimulation}
In this section,   two numerical simulations are presented to demonstrate  the
effectiveness of the theorems established in the previous section.

\begin{figure}[!htp]
\centering
 \includegraphics[height=1.5in]{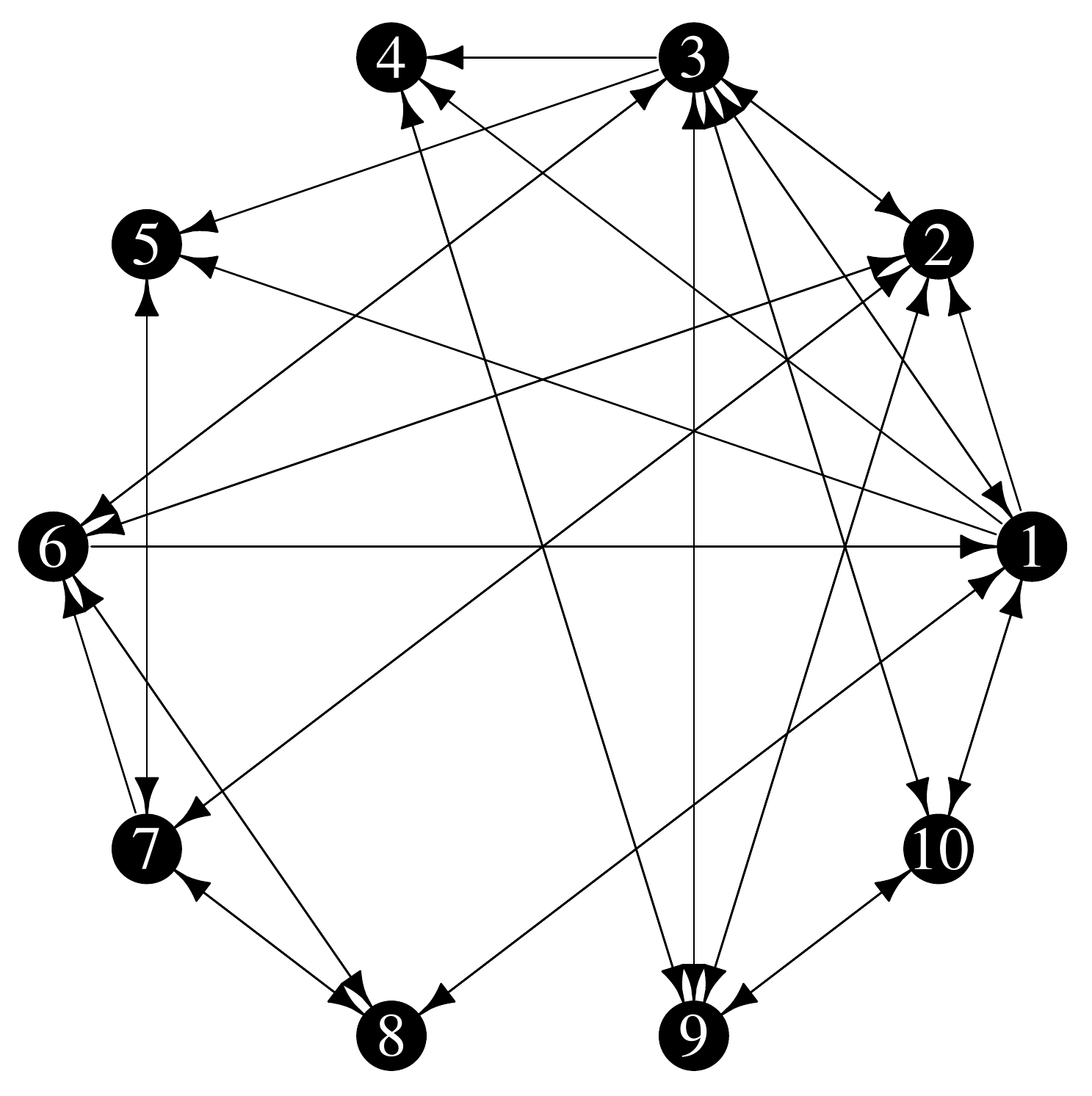}
 \caption{Directed graph $\gG$, where $a_{ij}=1$ for all $(j,i)\in \gE$ and $a_{ij} =0$ otherwise.}\label{fig:topstr}
\end{figure}

Consider a directed network of $10$ agents moving in the one-dimensional Euclidean space with the topology $\gG$
as shown in Fig. \ref{fig:topstr}. 

\begin{figure*}[!htp]
 \center
 \subfigure[$r_i(t)$]{\label{fig1:rit}\includegraphics[width=2.2in]{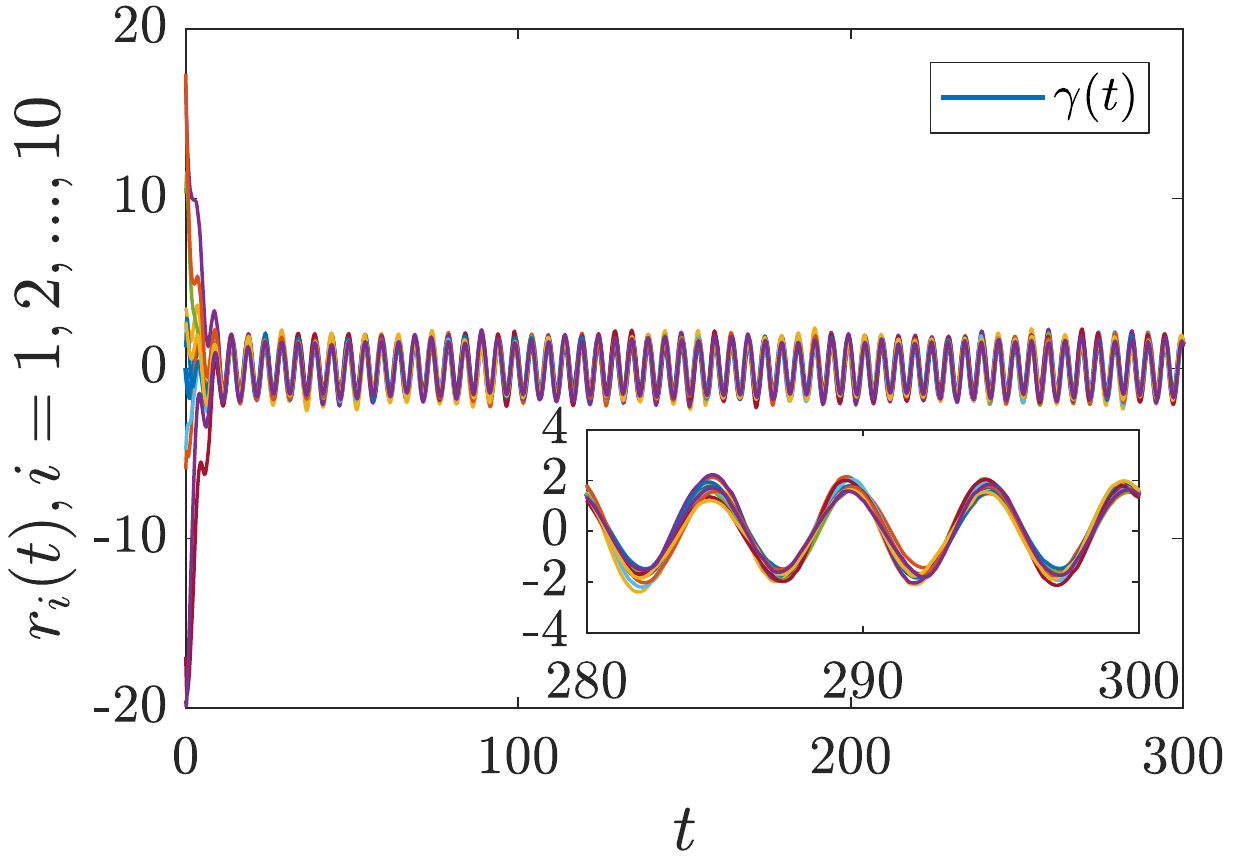}}~~~
 \subfigure[$v_i(t)$]{\label{fig1:vit}\includegraphics[width=2.2in]{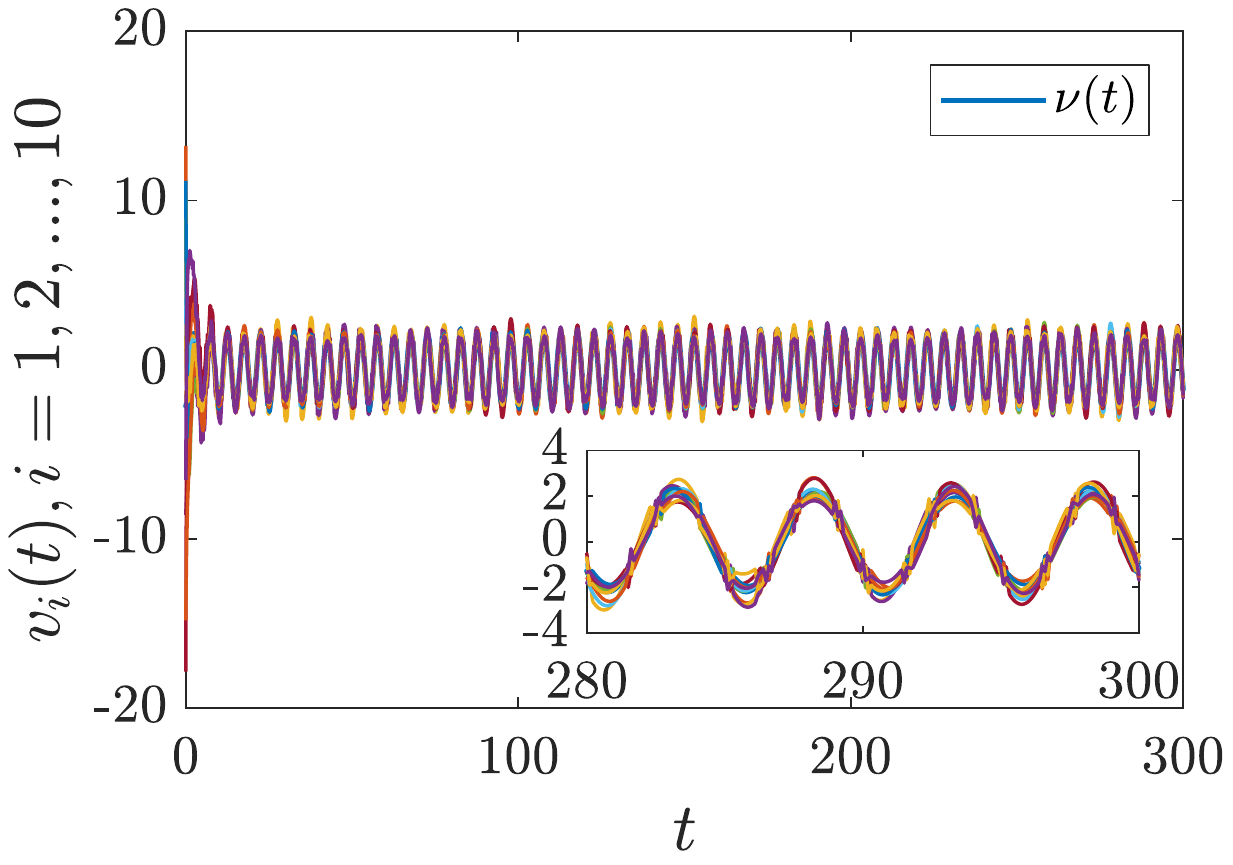}}~~~
 \subfigure[$E(t)$]{\label{fig1:et}\includegraphics[width=2.2in]{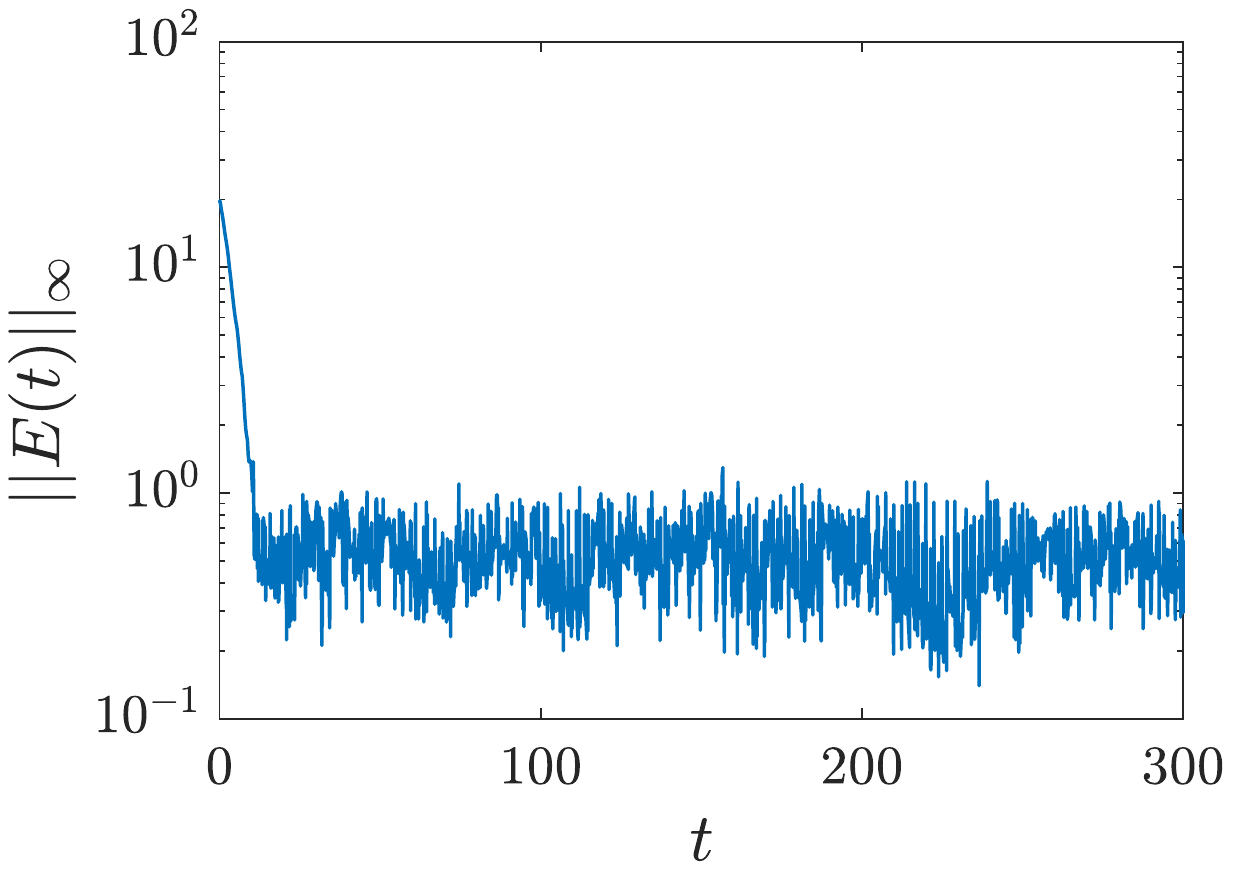}}
 \caption{Time evolutions of $r_i(t)$, $v_i(t)$ and $||E(t)||_{\infty}$ with a fixed zooming parameter $\mu$.}
 \label{fig1}
\end{figure*}

First,  consider the system under protocol \eqref{equ:quantcont} with a fixed zooming parameter (Case 1).
Let $\omega = \sqrt{\pi/2}$, thus $\tau \in ({k\pi}/{\omega},{k\pi}/{\omega}+\tau_0)$ for all $k = 0,1,2,\ldots$, where $\tau_0=\min{\textrm{arccot}(\phi_i)}$.
The relationship between $\lambda_i$ and $\phi_i$ is shown in Table \ref{tab1alphatau}.
Choose $\tau=0.1$, $M=10$ and $\Delta=0.5$.
Direct computation yields $\spec{\hat{\mE}} = 0.9747$, so that the conditions of Theorem \ref{thm1} are satisfied.
Define the synchronization error by $E(t) = [r(t)^{\top},v(t)^{\top}]^{\top} -  [\gamma(t),\nu(t)]^{\top} \otimes \mOne_n$.
The evolutions of the oscillators are shown in Fig. \ref{fig1} for the chosen initial values $(r_i(0),v_i(0))$, $i=1,2,\ldots,10$.
One can observe from Fig. \ref{fig1} that  complete synchronization cannot been reached, although the states converge to a bounded region of the synchronized orbits.

%

\begin{figure*}[!htp]
 \center
 \subfigure[$r_i(t)$]{\label{fig3:rit}\includegraphics[width=2.2in]{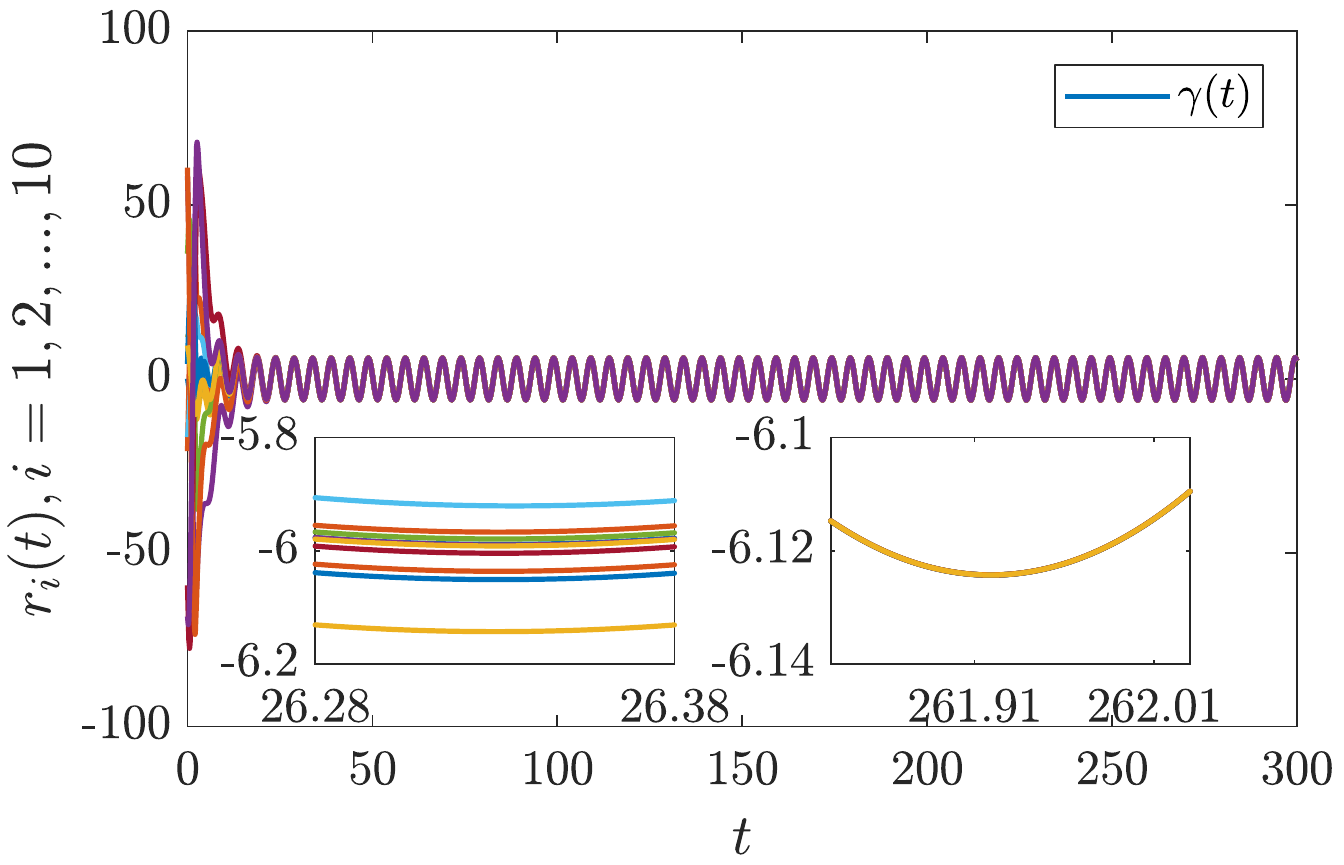}} ~~
 \subfigure[$v_i(t)$]{\label{fig3:vit}\includegraphics[width=2.2in]{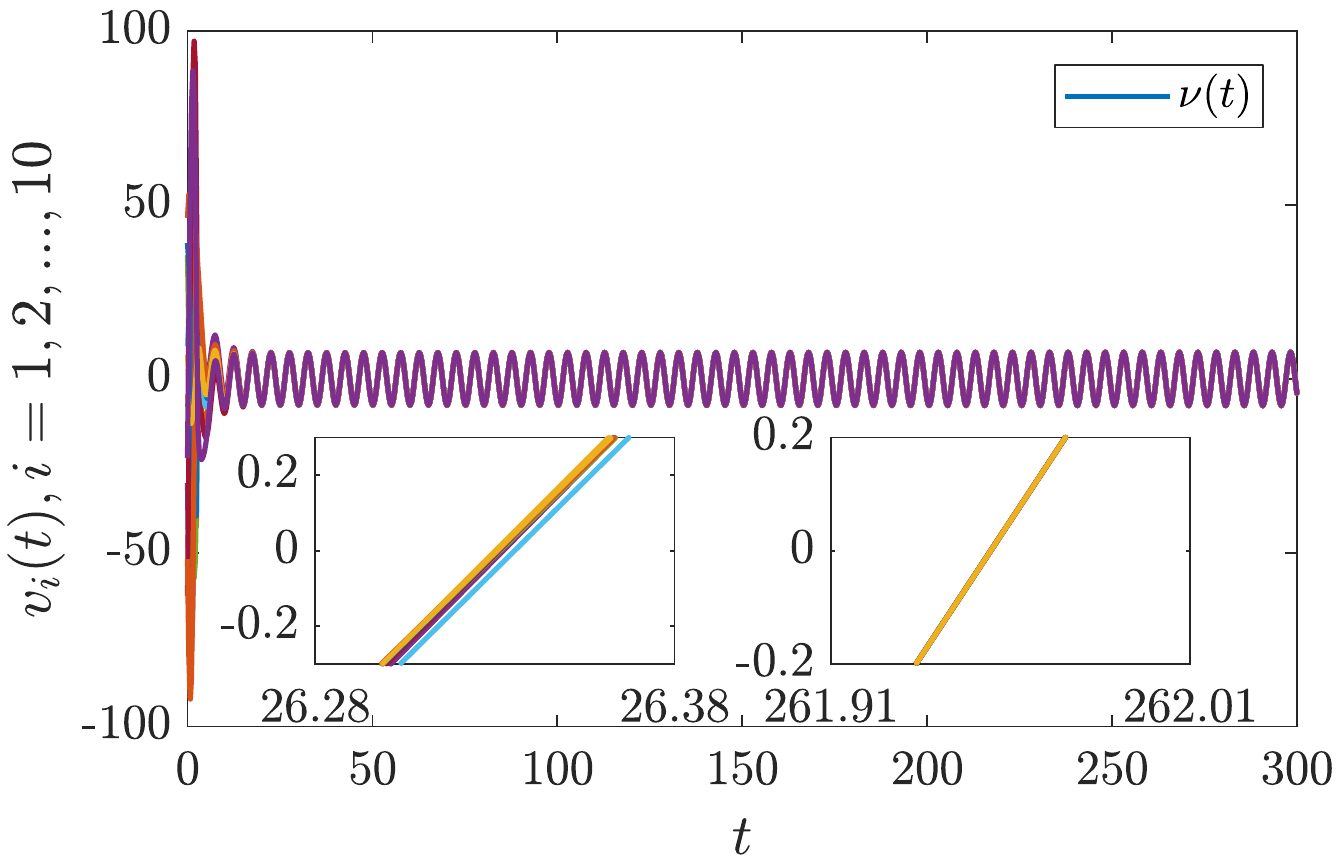}} ~~
 \subfigure[$\norminf{E(t)},\mu(t)$]{\label{fig3:et}\includegraphics[width=2.2in]{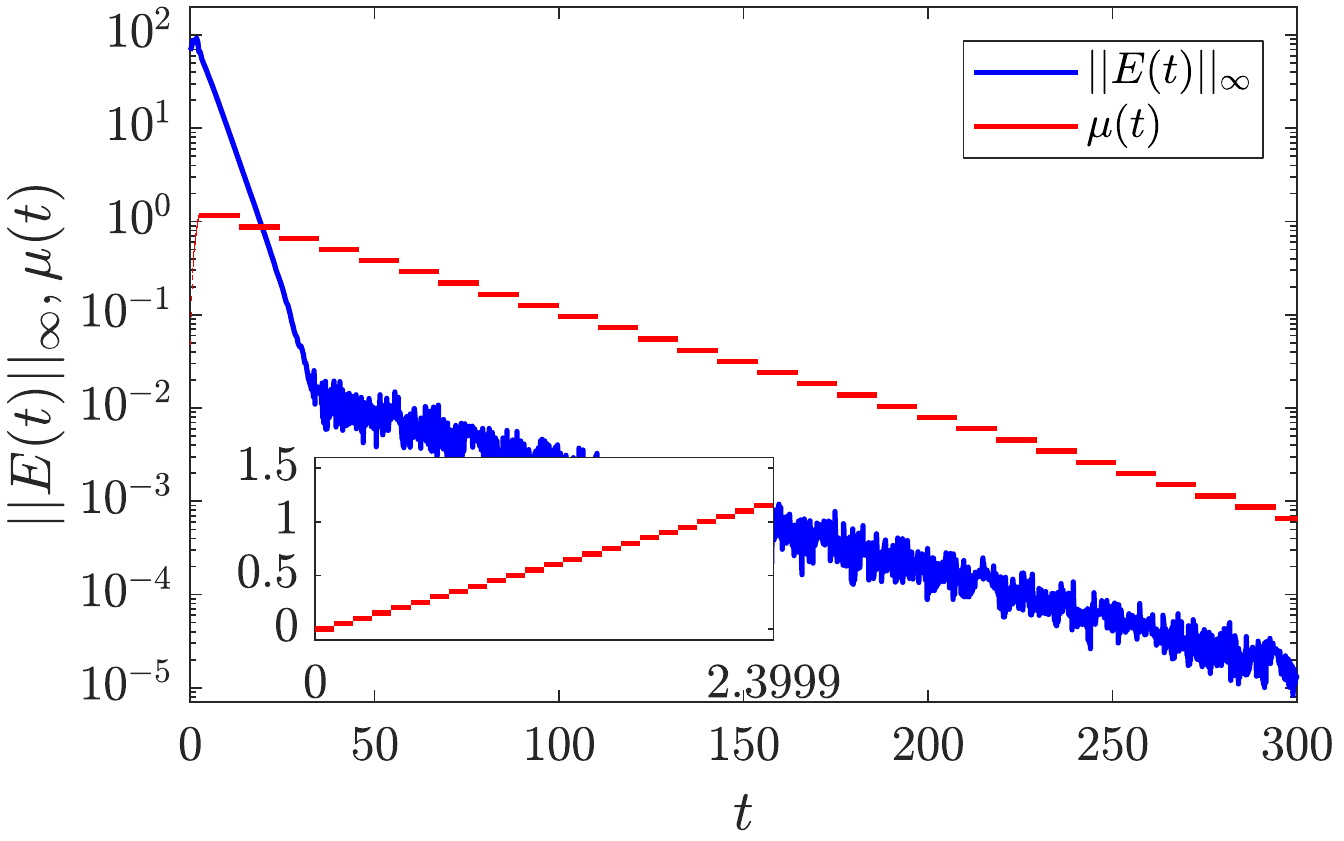}} 
 \caption{Time evolutions of $r_i(t)$, $v_i(t)$ and $||E(t)||_{\infty}$ with an adjustable zooming parameter $\mu(t)$.}
 \label{fig3}
\end{figure*}

\begin{table}[!htp]
\center
\caption{computation results.}
\begin{tabular}{ccccc}
 \hline
 $i$ & $\lambda_i$ & $\phi_i$ & $\text{arccot}(\phi_i)$\\
 \hline\hline
 2& $1.5594$ &1.2442 & 0.6770\\
 3& $6.3182 + 0.0706\mi$ &5.0419&0.1958\\
 4& $6.3182 - 0.0706\mi$ &5.0419&0.1958\\
 5& $2.9473$ &2.3516 &0.4021\\
 6& $3.4893 + 0.2867\mi$ &2.8052&0.3424\\
 7& $3.4893 - 0.2867\mi$ &2.8052&0.3424\\
 8& $5.1342$ &4.0965&0.2394\\
 9& $4.7440$ &3.7852&0.2583\\
 10& $3.0000$ &2.3937&0.3957\\
 \hline
\end{tabular}
\label{tab1alphatau}
\end{table}

Next, consider the system under protocol \eqref{equ:quantcont} with an adjustable zooming parameter (Case 2), where $u_i(t)$ and $\mu(t)$ are defined in equation \eqref{equ:uimuconc}.
Fig. \ref{fig3} shows the evolutions of the oscillators for the same initial values.
It is clear from Fig. \ref{fig3} that the complete synchronization is achieved.

\section{Conclusion}\label{sec:conclusion}
In this technical note, an effective quantized sampled-data feedback coupling protocol has been designed and evaluated for synchronizing networked harmonic oscillators.
The quantizer with a fixed or an adjustable zooming parameter has also been designed by using only sampled velocity data.
Some sufficient conditions have  been established under which the networked harmonic oscillators could achieve complete synchronization.

Future studies may include the synchronization of some general complex dynamical systems via quantized control, and 
the synchronization by designing coupling and control protocols with logarithmic quantizers.



\end{document}